  \documentclass{amsart}

   \usepackage{rotating}
  \usepackage{floatpag}
  \rotfloatpagestyle{empty}

  \usepackage{amsthm}
  \usepackage{amssymb, amsmath}

  \usepackage{graphicx}
    \usepackage{tikz}
  \usepackage{multicol}

  \theoremstyle{plain}
  \newtheorem{theorem}{Theorem}[section]
  \newtheorem{lemma}[theorem]{Lemma}
  \newtheorem{corollary}[theorem]{Corollary}

  \theoremstyle{definition}
  
  \newtheorem{example}[theorem]{Example}
\newtheorem{assumption}[theorem]{Assumption}

  \theoremstyle{remark}
  \newtheorem*{remark}{Remark}
  
 \newtheorem{question}{Question}
  \newtheorem{conjecture}{Conjecture}
  
\newcommand{\codim}{\operatorname{codim}}

\newcommand{\Gv}{G_{\hat{v}}}

\newcommand{\PP}{\ensuremath{\mathbb{P}}}

\newcommand{\RR}{\ensuremath{\mathbb{R}}}

\newcommand{\mI}{\ensuremath{\mathcal{I}}}

\begin{document}

\title[Line arrangements]{Line arrangements modeling curves of high degree: equations, syzygies and secants}
\author{Gregory Burnham, Zvi Rosen, Jessica Sidman, Peter Vermeire}

\email{ghburnham@gmail.com}

\address{Department of Mathematics, University of California, Berkeley, CA 94720}
\email{zhrosen@math.berkeley.edu}

\address{Department of Mathematics and Statistics, Mount Holyoke College, South Hadley, MA 01075}
\email{jsidman@mtholyoke.edu}

\address{Department of Mathematics, Central Michigan University, Mount Pleasant, MI, 48859}
\email{p.vermeire@cmich.edu}

\begin{abstract}
We study curves consisting of unions of projective lines whose intersections are given by graphs.  Under suitable hypotheses on the graph, these so-called \emph{graph curves}  can be embedded in projective space as line arrangements.   We discuss property $N_p$ for these embeddings and are able to produce products of linear forms that generate the ideal in certain cases.   We also briefly discuss questions regarding the higher-dimensional subspace arrangements obtained by taking the secant varieties of graph curves.
\end{abstract}

\maketitle


\section{Introduction}
An arrangement of linear subspaces, or subspace arrangement, is the union of a finite collection of linear subspaces of projective space.   In this paper we study arrangements of lines called graph curves with high degree  relative to genus.   We are particularly interested in the defining equations and syzygies of these subspace arrangements.  We will assume an algebraically closed ground field of characteristic zero throughout.

Let $G = (V, E)$ be a simple, connected graph with vertex set $V$ and edge set $E$.   Following \cite{cilibertoMiranda}, we assume that $G$ is \emph{subtrivalent}, meaning that each vertex has degree at most three. The (abstract) \emph{graph curve} $C_G$ associated to $G$ is constructed by taking the union of $\{L_v \mid v \in V\}$ where each $L_v$ is a copy of $\PP^1$ and lines $L_u$ and $L_v$ intersect in a node if and only if there is an edge between $u$ and $v$ in $G$.   (Note that if we think of the nodes of $C_G$ as vertices and the lines $L_v$ as edges, then $C_G$ is the graph dual to $G$.)   Since we are assuming that each vertex has degree less than or equal to three, $C_G$ is specified by purely combinatorial data; we may assume that on each component of $C_G$ the nodes are at 0, 1 or $\infty$.  Note that if each vertex of $G$ is trivalent, then each copy of $\PP^1$ in $C_G$ contains three nodes, and $C_G$ is stable.  (See \cite{bayerEisenbud, cilibertoMiranda}.)  

The motivation for the work presented here was to see if the syzygies of a high degree graph curve and its secant varieties would behave as they are expected to when the curve is smooth.  The $k$th secant variety, $\Sigma_k$, of a smooth curve in $\PP^r$ has expected dimension $\min \{2k+1, r\}$.  Thus, we expect the $k$th secant variety of $C_G$ to be an arrangement of subspaces of dimension $2k+1.$ 

Many authors (\cite{einLazarsfeld2, gallegoPurnaprajna, heringSchenckSmith, ottavianiPaoletti}) have given generalizations of the results for smooth curves to higher-dimensional varieties, showing that embeddings via line bundles satisfying various positivity conditions will also satisfy property $N_p$.  However, recent work of Ein and Lazarsfeld \cite{einLazarsfeld2} shows that these results describe only a small portion of the minimal free resolution of a higher-dimensional variety, and what happens in the remaining piece is quite complicated, contrary to the belief that positivity of an embedding simplifies syzygies.

One can view the conjecture of \cite{vermeire}, which says that we should expect property $N_{k+2,p}$ (ideal generators of degree $k+2$ and linear syzygies through stage $p$) for the $k$th secant variety of a smooth curve of genus $g$ embedded via a complete linear series of degree at least $2g+2k+1+p$ as an alternate way of generalizing property $N_p$ for curves to higher-dimensional varieties.  Some progress was made for first secant varieties, using geometric methods in \cite{sidmanVermeire}, but the recursive nature of these methods makes generalizing those techniques to higher secant varieties daunting.  If a similar result were true for the secant varieties of graph curves, the proof methods would necessarily be very different and the hope is that they would shed new insights into understanding secant varieties of smooth curves.

Based on many examples computed with \emph{Macaulay2} \cite{gs}, the situation looks promising.  However, when $g > 2$, the combinatorics can be intricate even if we only consider curves and not secant varieties.  We will focus on curves in \S \ref{Np} and \ref{explicitEmbed} and turn to a discussion of the syzygies of secant varieties of graph curves in \S \ref{secants}.

We begin by setting some assumptions and notation.  Let $d$ be the number of vertices in $G$. The topology of $G$ determines the arithmetic genus of $C_G$ as we may view $G$ as a 1-dimensional simplicial complex, from which it follows that $p_a(C_G) = h^1(G, k)$ if $G$ is connected (see Proposition 1.1 in \cite{bayerEisenbud}).  We refer to this quantity as the genus $g$ of $G$, and $|E|=d+g-1$.  Note that $g$ is not the genus of $G$ in the usual graph-theoretic sense. 

The story that we wish to generalize to the setting of graph curves began with Green and Lazarsfeld \cite{greenLazarsfeld} in the early 80's who showed that if $C$ is a smooth and irreducible curve of genus $g$ embedded in projective space via a complete linear series of degree $d \geq 2g+1+p$, then $C$ satisfies property $N_p.$  In other words, its ideal is generate by quadrics with syzygy modules generated by linear forms through the $p$th stage of the resolution. 

 We conjecture that if $G$ satisfies Assumption \ref{recursive}, then 
property $N_p$ will hold for $C_G$ embedded as a line arrangement in $\PP^{d-g}.$

\begin{assumption}\label{recursive}
Fix $p \geq 0,$ and let $G$ be a simple, connected, subtrivalent graph with $d \geq 2g+1+p$.  Assume that if $G'$ is a connected subgraph induced on $V' \subset V,$ $d' = |V'|,$ and $g'$ is the genus of $G',$ then $d' \geq 2g'+1+p$ if $g' \geq 1$
\end{assumption}

To see that the recursive hypotheses are necessary, note that a graph may satisfy $d \geq 2g+2$, but if it contains a triangle, then the ideal of $C_G$ cannot be generated by quadrics.

If Assumption \ref{recursive} is satisfied for some $p \geq 0,$ then $C_G$ embeds in $\PP^{d-g}$ as a line arrangement via \cite{cfhr} and is arithmeticaly Cohen-Macaulay by \cite{franciosiTenni}.  If $C_G$ is arithmetically Cohen-Macaulay, we may proceed as in \cite{greenLazarsfeld} and property $N_p$ for $C_G$ will follow if property $N_p$ holds for a general hyperplane section.   In \cite{greenLazarsfeld}, Green and Lazarsfeld deduce property $N_p$ for points in linearly general position, and conjectured that the failure of property $N_p$ for a set of $2r+1+p$ points implied the existence of a subset of $2k+2-p$ points on a $\PP^k.$  As shown in \cite{eisenbudKoh, green}, this conjecture for point sets is a consequence of the linear syzygy conjecture of Eisenbud, Koh, and Stillman \cite{eks}.  Green proved the linear syzygy conjecture in \cite{green}, and for graph curves of degree $g \leq 2$ we can show that an embedding of $C_G$ as a line arrangement via a complete linear series must satisfy $N_p$ if Assumption \ref{recursive} is satisfied.

%

Graph curves associated to graphs in which every vertex is trivalent are canonical curves, and have been studied in several different contexts.  For example, Ciliberto, Harris, and Miranda \cite{cilibertoHarrisMiranda} used graph curves to understand the surjectivity of the Wahl map, Ciliberto and Miranda \cite{cilibertoMiranda} related graph curves to graph colorings, and Bayer and Eisenbud \cite{bayerEisenbud} studied graph curves in connection with Green's conjecture.  In fact, Proposition 3.1 in \cite{bayerEisenbud} gives an explicit description of generators of the ideal of a canonical graph curve using the combinatorics of $G$.  More recently Ballico has written several papers about graph curves \cite{ballico1, ballico2}.  

We present an explicit embedding of $C_G$ into projective space in \S \ref{explicitEmbed}.  If the ideal of $C_G$ is generated by quadrics, this allows us to show that $I_{C_G}$ may be generated by products of linear forms (Theorem~\ref{thm: products}).   

Although a subspace arrangement may always be cut out by products of linear forms set-theoretically, we do not generally expect the ideal of a subspace arrangement to be generated by products of linear forms, cf Proposition 5.4 and Proposition 5.7 in \cite{bps}.  
The most interesting examples of subspace arrangements with ideals generated by products of linear forms occur when the intersections among the subspaces have a rich combinatorial structure \cite{lili, lovasz,bps}.  
If $G$ is a path or a cycle, then $C_G$ can be embedded in projective space so that its ideal is generated by square free monomials.  In both cases, the ideals of the nontrivial secant varieties of these curves are also generated by square free monomials and are examples of ``combinatorial secant varieties" \cite{sturmfelsSullivant}. 

In addition to viewing graph curves and their secant varieties as combinatorial models of smooth curves and their secant varieties, we can also 
think of them as a new way of generating arrangements of linear subspaces with interesting interactions between the combinatorics of the arrangements, the geometry of the embeddings, and their defining equations.  We present conjectures and questions for further work in this direction in \S \ref{questions}.
\bigskip

\noindent{\bf Acknowledgements}
The third author thanks the Abel Symposium ``Combinatorial aspects of commutative algebra and algebraic geometry" from 2009, where David Eisenbud and Mike Stillman suggested the problem of considering graph curves and their secant varieties.  We thank Ian Barnett, David Cox, David Eisenbud, Marco Franciosi, Jenia Tevelev for helpful correspondences.  We are very grateful to an anonymous referee who led us to make significant and necessary changes to the original version of this paper\footnote{Readers of \cite{sidmanSchenck} should note that references there referring to this paper refer to the original version, and results are not in the same place or form here.}.   This work began at the Mount Holyoke College Summer Mathematics Institute funded by NSF grant DMS-0849637.  The third author was supported by NSF grant DMS-0600471.  Finally the third author wishes to thank Rob Lazarsfeld for introducing her to algebraic geometry, syzygies and subspace arrangements.

\section{Regularity and  property $N_p$}
\label{Np}

In this section we will show that if $g \leq 2,$ then the ideal of a linearly normal embedding of $C_G$ as a line arrangement satisfies property $N_p$ if $G$ satisfies Assumption \ref{recursive} for some $p \geq 0$ following the idea of the ``quick" proof that a smooth and irreducible curve of degree $d \geq 2g+1+p$ satisfies $N_p$ given in \cite{greenLazarsfeld}.  

A key assumption in \cite{greenLazarsfeld} is that a hyperplane section of a smooth curve of degree $d \geq 2g+1+p$ will consist of points in linearly general position.  This fact is used to show that the points in a hyperplane section of the curve impose independent conditions on quadrics.  

Using Lemma \ref{cycle} we can show that this is not the case for a graph curve if $G$ contains a cycle as a proper subgraph.  

\begin{lemma} \label{cycle}
If $G$ is a cycle on $d$ vertices, then a hyperplane section of $C_G$ has a 1-dimensional space of linear dependence relations and all of the points are contained in the support of the relation.
\end{lemma}
\begin{proof}
A cycle of length $d$ embeds into $\PP^{d-1}$, so the hyperplane section consists of $d$ points in $\PP^{d-2}.$  A set of $d$ points spanning a $\PP^{d-2}$ must satisfy a unique relation up to scalar.
\end{proof}

Therefore, if we have a cycle as a proper subset of a graph $G$, the points of a hyperplane section must fail to be in linearly general position.  Because $N_p$ fails if $G$ contains a cycle of length $p+2$, it will often be impossible to reproduce the graded Betti diagrams of a smooth curve with the graded Betti diagrams of a graph curve.  For instance, for genus $g=2$ and degree $d$, the length of the smallest cycle has an upper bound of $\lfloor\frac{2d-1}{3}\rfloor + 1$. 

Nevertheless, we will show that if $G$ satisfies Assumption \ref{recursive}  for $g\leq 2$, then a general hyperplane section of $C_G$ imposes independent conditions on quadrics.  This follows from the weaker assumption that no $2k+2$ of the points lie on a $\PP^k$ using ideas from \cite{eisenbudKoh}.

\begin{theorem}\label{thm: points}
Suppose that $G$ satisfies Assumption \ref{recursive} for some $p \geq 0$, $g \leq 2$, and $C_G$ is embedded in $\PP^{d-g}$ as a line arrangement via a complete linear series.  If H is a general hyperplane and $X = H \cap C_G,$ then there is no set of $2k+2-p$ dependent points of $X$ lying on a $\PP^k.$
\end{theorem}
\begin{proof}
Let $Y \subset X$. Suppose for contradiction that $|Y| = 2k+2-p$ and $Y$ spans a $\PP^k.$  This means that there is a $2k+2-p-(k+1) = k+1-p =m$ dimensional space of dependence relations on $Y.$  Since $g\leq 2,$ we know that  $m \leq 2.$  If $m = 0,$ then the points are independent which contradicts our hypotheses.  


If $m = 1,$ then $k = p.$  Either the support of the unique dependence relation on $Y$ contains a cycle of points, or the relation is a linear combination of dependence relations on 2 cycles in which at least one point has been eliminated from their support.   If $\{ \gamma_i\}$ form a basis for $H_1(G;\RR),$ the corresponding dependence relations $\{R_i\}$ form a basis for the space of linear relations on $X$, and Assumption \ref{recursive} implies that $\gamma_1 \cup \gamma_2$ contains at least $5+p$ points.    The cycles $\gamma_1$ and $\gamma_2$ can be combined in $H_1(G; \RR)$ to form a distinct cycle $\gamma_3$ which also supports a unique linear dependence. Therefore, if we fix the coefficient of $R_1$ there is a unique multiple of $R_2$ that eliminates the shared points in the interior of their common path to create a dependence relation with support on $\gamma_3.$  Consequently, we see that we cannot simultaneously eliminate the endpoints of this path and the points between them from the support.  Therefore, if a linear combination of $R_1$ and $R_2$ is not supported on a full cycle, it contains at least $2\cdot 2+1+p = 2+1+p = 3+p$ points, implying that $Y$ spans a projective space of dimension at least $p+1,$ which is a contradiction as $k=p.$

If $m = 2,$ then $k = p+1.$  In this case $g=2$, and $Y$ must contain the support of both cycles of $G$, in which case $2k+2-p \geq 2\cdot 2+1+p,$ or $2k \geq 2p +3,$ which contradicts $k=p+1.$
\end{proof}

\begin{remark}
We conjecture that if $G$ satisfies Assumption \ref{recursive} and $C_G$ is embedded via a complete linear series then Theorem \ref{thm: points} holds for all $g.$  The idea is that if there is an $m$-dimensional space of dependence relations on $Y$, then we need at least $m$ independent cycles of $G$ to span this space.  The support of $m$ cycles contains at least $2m+1+p$ points.   If more than $m$ cycles are needed to span the space of dependence relations of $Y$, then we may have eliminated some points from the support, but we will always have at least $2m+1+p$ points remaining.  
\end{remark}

\begin{theorem}\label{points}
If $G$ satisfies Assumption \ref{recursive}  for some $p\geq 0$,  and no $2k+2-p$ points of $X$ lie on a $\PP^k,$  then a general hyperplane section of $C_G$ has a 3-regular ideal and satisfies property $N_p$.
\end{theorem}
\begin{proof}
The proposition on pg. 169 of \cite{eisenbudKoh} states that $X$ imposes independent conditions on quadrics if $X$ does not contain a subset of $2k+2$ points on a projective $k$-plane.  This implies that the ideal of $X$ is 3-regular by Lemma 2 of \cite{eisenbudKoh}.  The ideal of $X$ satisfies $N_p$ as a consequence of Theorem 2.1 in \cite{green}.
\end{proof}

\begin{theorem}\label{ACM}
Suppose that $G$ satisfies Assumption \ref{recursive} for some $p\geq 0$, $d \geq 2g+1+p,$ and $C_G$ is embedded in $\PP^{d-g}$ as a line arrangement via a complete linear series.  If no $2k+2-p$ points of a general hyperplane section lie on a $\PP^k,$ then this embedding is arithmetically Cohen-Macualay, 3-regular and satisfies $N_p$.
\end{theorem}

\begin{proof}
For $3$-regularity we need $H^1(\mI_{C_G}(2))=H^2(\mI_{C_G}(1))=0.$  We know that $H^1(O_{C_G}(1))=0$ by Serre duality and our hypothesis that $d\geq 2g+2$.  This implies that  $H^2(\mI_{C_G}(1))=0$.  To see the vanishing of $H^1(\mI_{C_G}(2))$, note  via Theorem \ref{points} the regularity of the ideal of a general hyperplane section $X$ of $C_G$ is 3 which implies that $H^1(\mI_X(2))=0.$  Since $C_G$ is embedded via a complete linear series, $H^1(\mI_{C_G}(1))=0,$ and we conclude that $H^1(\mI_{C_G}(2))=0.$

The curve $C_G \subset \PP^{d-g}$ is arithmetically Cohen-Macaulay if its homogeneous coordinate ring is Cohen-Macaulay.  Equivalently, the hypersurfaces of degree $m$ are a complete linear series, which holds if and only if $H^1(\mI_{C_G}(m))=0$ for all $m \geq 0.$  (See Section 8A of \cite{eisenbud}.)  When $m=0$, this follows because $C_G$ is connected.  We know that $H^1(\mI_{C_G}(1))=0$ from  the linear normality of the embedding and $H^1(\mI_{C_G}(k))=0$ for all $k \geq 2$ by the 3-regularity of the ideal.
\end{proof}

\begin{corollary}\label{gleq2}
If $G$ satisfies Assumption \ref{recursive} for some $p\geq 0$, and $g \leq 2,$ then an embedding of $C_G$ as a line arrangement via a complete linear series is arithmetically Cohen-Macualay, 3-regular and satisfies $N_p$.
\end{corollary}
\begin{proof}
Theorem \ref{thm: points} implies that the hypotheses of Theorems \ref{points} and \ref{ACM} are satisfied.
\end{proof}

Note that by Theorem 4.2 of \cite{franciosiTenni}, we know that if Assumption \ref{recursive} holds for some $p \geq 0$, then an embedding of $C_G$ as a line arrangement is always Cohen-Macualy, as our singularities are planar.
Moreover, Ballico and Franciosi \cite{ballicoFranciosi} proved that a line bundle $L$ on a reduced curve $C$ satisfies property $N_p$ under certain numerical conditions on the positivity of $L$ with respect to subcurves constructed from an ordering of the irreducible components of $C.$  Their hypothesis on the degree of $L$ restricted to an irreducible component fails if $G$ contains a cycle or if $p >0$, and $L$ has degree 1 on each line.   However, if $G$ is a tree, then Assumption \ref{recursive} is automatically satisfied, so we expect that the ideal of $C_G$ is 2-regular in this case.  In fact, this follows from \cite{eghp} because the lines in $C_G$ can be ordered in such a way that the $i$th line intersects the span of the previous lines in a single point.

\section{Line arrangements generated by products of linear forms} \label{explicitEmbed}

In this section we present an embedding of $C_G$ into projective space if its edges can be labeled according to certain rules described below.  If the ideal of $C_G$ is generated by quadrics, then we identify conditions on the labeling that guarantee the existence of generators of the ideal of $C_G$ that are monomial and binomial products of linear forms.


Given a graph $G$ satisfying Assumption \ref{recursive}, construct $\tilde{G}$ from $G$ by adding a loop to each vertex of degree 1 so that vertices of degree 1 in $G$ are incident to two edges in $\tilde{G}.$  For the induction in Theorem \ref{thm: products} we also need to allow the possibility of the addition of a loop at vertices with degree 2 in $G.$  We describe the embedding of $C_G \subset \PP^{d-g}$ by labeling the edges of $\tilde{G}$ with monomial and binomial linear forms in $S[x_0, \ldots, x_{d-g}]$ that indicate how coordinates of $\PP^{d-g}$ parameterize each line $L_v.$  

Label each edge of $\tilde{G}$ with a monomial $x_i$ or a binomial $x_i-x_j$ subject to the following rules:
\begin{enumerate}
\item   We require that each variable $x_i$ appears as a monomial edge label exactly once.
\item  Binomials only appear on non loop edges.
\item  Each edge labeled with a binomial is incident to a vertex with 3 incident edges.
\item  If $v$ has 3 incident edges, then they are labeled $x_j$,$x_k$, and $x_j-x_k$, where $j\neq k \in  \{0, \ldots, d-g\}.$ 

\end{enumerate}
For a fixed graph $G$, it may be the case that some $\tilde{G}$ can be labeled according to these rules and others may not.

To define the ideal of $L_v$ let $\Omega_v$ be the set defined by deleting all of the variables appearing on the edges incident to $v$ from the set of variables of $S$ and then adding in the binomial edge label incident to $v$ if $v$ has only 2 incident edges in $\tilde{G}$.   We let $I_v =  \langle \Omega_v \rangle$ be the ideal of $L_v.$  Thus, the line $L_v$ is parameterized by the coordinates on the incident edges, with coordinates $i$ and $j$ equal if $x_i-x_j$ appears at $v$ but $x_i$ and $x_j$ do not.

\begin{example}\label{g2d5}
The graph $G$ below has $g=2$ and $d=5.$    
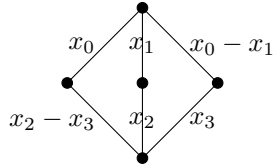
\begin{multicols}{2}

\[
	\begin{tikzpicture}
	\filldraw [black] (-1,0) circle (2pt);
	\filldraw [black] (0,0) circle (2pt);
	\filldraw [black] (1,0) circle (2pt);
	\filldraw [black] (0,1) circle (2pt);
	\filldraw [black] (0,-1) circle (2pt);

	\draw (-1,0)--( -.5,.5)node[left]{$x_0$} -- (0,1) --(0,.5)node[]{$x_1$}-- (0,0)-- (0,-.5)node[]{$x_2$}--(0,-1) -- (-.5,-.5)node[left]{$x_2-x_3$}--cycle;
	\draw (0,1) --(.5,.5)node[right]{$x_0-x_1$}-- (1,0)--(.5,-.5)node[right]{$x_3$} -- (0,-1);
	\end{tikzpicture}
\]

The ideals of the 5 lines are
\[
\begin{split}
 \langle  x_2, x_3\rangle \\
\langle  x_1, x_2-x_3 \rangle \\
\langle x_0,x_3 \rangle\\
\langle x_0-x_1, x_2\rangle\\
\langle x_0,x_1\rangle
\end{split}
\]

\end{multicols}
Via \emph{Macaulay2} \cite{gs}, the ideal of the arrangement is \[\langle {x}_{0} {x}_{2}-{x}_{0} {x}_{3}+{x}_{1} {x}_{3},{x}_{1} {x}_{2} {x}_{3},{x}_{0} {x}_{1} {x}_{3}-{x}_{1}^{2} {x}_{3} \rangle.\]
The labeling gives rise to an embedding of $C_G$, but the ideal of this embedding is not generated by products of linear forms and is not generated by quadrics.

\end{example}
\begin{theorem}
If $G$ satisfies Assumption \ref{recursive} for $p\geq 0$ and $\tilde{G}$ is labeled as described above, then $I = \cap_{v \in V} I_v$ is the ideal of an embedding of $C_G$ into $\PP^{d-g}.$
\end{theorem}
\begin{proof}
If $u$ and $v$ are connected by an edge in $G$ and $\ell$ is the linear form on the edge that joins them, then the lines $L_u$ and $L_v$ intersect at the point of $\PP^{d-g}$ that has the coordinates appearing in $\ell$ set to 1 and all other coordinates set to 0.  

To see that a labeling defines an embedding of $C_G$ we must show that if $u$ and $v$ are not connected by an edge, then $L_u$ and $L_v$ do not intersect.  If they intersect then a variable appearing on an edge incident to $u$ must also appear on an edge incident to $v.$  If $x_i$ is an edge label at $u$, and $v$ is incident to an edge with a label containing $x_i$, we must have the configuration on the left in Figure \ref{embedCase}.   But then the only coordinates of $L_v$ with $x_i$ nonzero also have $x_j$ nonzero, and $L_u$ does not contain any points with $x_j$ nonzero unless $w$ is trivalent and there is an edge labeled $x_j-x_k$ incident to $u$.  However, this is not possible, because $G$ does not contain any triangles.  Hence, the lines cannot intersect. 

\begin{figure}[htbp]
\begin{multicols}{2}

	\begin{tikzpicture}
		\filldraw [black] (-1,0) circle (2pt);
		\filldraw [black] (1,0) circle (2pt);
		\filldraw [black] (0,1) circle (2pt);
		\filldraw [black] (0,2) circle (2pt);

		\draw (-1,0)node[left]{$u$} --(-.5, .5)node[left]{$x_i$}--(0,1)--(.5,.5)node[right]{$x_i-x_j$} -- (1,0)node[right]{$v$};
		\draw (0,1) --(0, 1.5)node[left]{$x_j$}-- (0,2)node[above]{$w$};
	\end{tikzpicture}
		
	\begin{tikzpicture}
		\filldraw [black] (-1,0) circle (2pt);
		\filldraw [black] (1,0) circle (2pt);
		\filldraw [black] (0,1) circle (2pt);
		\filldraw [black] (0,-1) circle (2pt);
		\filldraw [black] (1,1) circle (2pt);
		\filldraw [black] (-1,-1) circle (2pt);

		\draw (-1,0) node[left]{$u$}--(-.5,.5)node[left]{$x_i-x_j$}-- (0,1) --(.5,1)node[above]{$x_j$}-- (1,1);
		\draw (0,1) -- (0,0)node[left]{$x_i$}--(0,-1)--(.5, -.5)node[right]{$x_i-x_k$} -- (1,0)node[right]{$v$};
		\draw (-1,-1)--(-.5, -1)node[below]{$x_k$} -- (0,-1);
	\end{tikzpicture}

\end{multicols}
\caption{}
\label{embedCase}
\end{figure}

The only other possibility is that $x_i$ appears in a binomial at $u$ and at $v$ as in the diagram on the right in Figure \ref{embedCase}.  But then if the lines $L_u$ and $L_v$ intersect, the three coordinates $x_i, x_j, x_k$ must all be nonzero and equal.  This means that the edge labeled with $x_k$ must be incident to $u$ and the edge labeled $x_j$ must be incident to $v.$  But this is forbidden because $d \geq 2g+1$ for all subgraphs of genus $g.$
\end{proof}


Our method of labeling edges with linear forms is similar in spirit to the description of the generators of the ideal of a canonical graph curve (corresponding to a trivalent graph) in \cite{bayerEisenbud}.  They label the edges in $G$ with a basis for the space of 1-cochains of $G$ and intersect an ideal generated by monomials in this basis with the ring generated by the 1-cocycles.

In order to describe the generators of $I_{C_G}$ explicitly, we must make some further assumptions on the relative placement of labels.

\begin{assumption}  \label{ass:label}

The labeling on $\tilde{G}$ satisfies the following conditions.
\begin{enumerate}

\item  Incident edges never both have binomial labels.  In other words. the labeling below never appears.


	
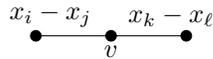
\begin{tikzpicture}
	\filldraw [black] (0,0) circle (2pt);
	\filldraw [black] (1,0) circle (2pt);
	\filldraw [black] (2,0) circle (2pt);

	\draw (0,0) -- (.2,0)node[above]{$x_i-x_j$} -- (1,0)node[below]{$v$} -- (1.8,0)node[above]{$x_k-x_{\ell}$} -- (2,0);
		\end{tikzpicture}

\item  If $v$ is a vertex of degree 2 as depicted below (with $i,j,k$ distinct), then there are no other edges with labels containing $x_i$ that are incident to edges with labels containing $x_j$ or $x_k.$ 

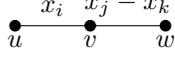
	\begin{tikzpicture}

	\filldraw [black] (0,0) circle (2pt);
	\filldraw [black] (1,0) circle (2pt);
	\filldraw [black] (2,0) circle (2pt);

	\draw (0,0)node[below]{$u$}-- (.5,0)node[above]{$x_i$} -- (1,0)node[below]{$v$} -- (1.5,0)node[above]{$x_j-x_k$} -- (2,0)node[below]{$w$};
	
	\end{tikzpicture}
\item  The vertices of $G$ are ordered $v_1, \ldots, v_d$, $G_i$ is the graph  induced on $v_1, \ldots, v_i$ and $\tilde{G}_{i-1}$ is obtained from $\tilde{G}_i$ by removing $v_i$ and  replacing any non loop edge $uv_i$ labeled with a monomial by a loop at $u$ labeled with the same monomial.
	\begin{enumerate}
	\item  $G_i$ is connected;
	\item  $v_i$ has at most 2 incident edges in $\tilde{G}_i$;
	\item  if $v_i$ is connected to a vertex $u$ in $G_{i-1}$ via an edge labeled with a binomial, then $u$ is incident to 3 in $G_i$.  (i.e, $L_u$ has a monomial ideal.)

	\end{enumerate} 
	\end{enumerate}
\end{assumption}
 
 In what follows, let $\Gv$ denote the subgraph of $G$ obtained by removing $v$ and all of its incident edges.   If $C_G$ is embedded in $\PP^{d-g},$ we let $C_{\Gv}$ be the corresponding subcurve.   Note that if $\deg v = 1$ in $G$ and we remove the line $L_v$ from $C_G$ embedded as above, then $C_{\Gv}$ is embedded as a line arrangement via a complete linear series in a hyperplane.  If $\deg v = 2$ in $G$ and $v$ is contained in a cycle, then $\Gv$ is still connected, the genus drops by 1, and the remaining subcurve is embedded via a complete linear series.  We do not allow the removal of vertices of degree 3 because if $\deg v = 3,$ and $\Gv$ is connected, then the genus drops by 2, and $C_{\Gv}$ is not embedded via a complete linear series.
 
 \begin{lemma}\label{gens}
 Suppose that Assumption \ref{recursive} holds for some $p \geq 1$ and Assumption \ref{ass:label} also holds.  If the configuration in part (2) of Assumption \ref{ass:label} appears in a labeling of $\tilde{G}$, then $x_i(x_j-x_k)$ is in the ideal of $C_G$ and $x_ix_j$, $x_ix_k$ are in the ideal of $C_{\Gv}$. \end{lemma}
\begin{proof}
To see that $x_i(x_j-x_k)$ is in the ideal of $C_G,$ note that $x_j-x_k$ is in the ideal of $L_v$.  Our hypotheses imply that for any vertex $v' \neq v$ that if $x_i$ appears on an incident edge, neither $x_j$ nor $x_k$ do, so $x_j-x_k$ is in the ideal of $L_{v'}.$  Otherwise, $x_i$ is in the ideal of $L_{v'}.$  Hence, $x_i(x_j-x_k)$ is in the ideal of each line.


It is easy to see that neither $x_ix_j$ nor $x_ix_k$ vanish on $L_v.$  Since this is the only line where the coordinate $x_i$ is paired with $x_j$ or $x_k$, it follows that these two monomials are contained in the ideal of $C_{\Gv}.$

  \end{proof}
 
 \begin{example}\label{g2}
If $g=2$, $G$ has precisely 2 trivalent vertices, and it satisfies Assumption \ref{recursive} for some $p \geq 1$, then it can be labeled according to Assumption \ref{ass:label}.   If the cycles are disjoint, then $G$ must consists of 2 cycles and a bridge between them.  Putting one binomial label in each cycle satisfies Assumption \ref{ass:label} because each cycle has length at least 4..

If the cycles overlap, then we have 3 paths between trivalent vertices $u$ and $v.$  Label the shortest path with monomials and put one binomial label on each of the remaining paths.   For example, the graph in Figure \ref{g2d6} satisfies Assumption \ref{ass:label}, and has defining ideal $\langle {x}_{3} {x}_{4},{x}_{0} {x}_{4}-{x}_{2} {x}_{4},{x}_{0} {x}_{3}-{x}_{1} {x}_{3},{x}_{1} {x}_{2}\rangle.$

\begin{figure}[htbp]
\begin{tikzpicture}
	\filldraw [black] (0,0) circle (2pt);
	\filldraw [black] (1,0) circle (2pt);
	\filldraw [black] (-1,0) circle (2pt);
	\filldraw [black] (1,1) circle (2pt);
	\filldraw [black] (-1,1) circle (2pt);
	\filldraw [black] (0,1) circle (2pt);

	\draw (-1,0) --(-.5,0)node[below]{$x_2$}--(.5,0)node[below]{$x_0-x_2$}-- (1,0) --(1,.5)node[right]{$x_4$}-- (1,1) -- (.5,1)node[above]{$x_1$}--(-.5,1)node[above]{$x_0-x_1$}--(-1,1) --(-1,.5)node[left]{$x_3$}-- cycle;
	\draw (0,0)--(0,.5)node[left]{$x_0$}--(0,1);
\end{tikzpicture}
\caption{}
\label{g2d6}
\end{figure}

\end{example}
%

\begin{theorem}\label{thm: products}
Suppose that $G$ satisfies Assumption \ref{recursive} with $p =1$.  Fix a $\tilde{G}$ and a labeling that gives an embedding of $C_G$ into $\PP^{d-g}$ as a line arrangement.  If Assumption \ref{ass:label} is satisfied,  and $I_{C_{G_i}}$ is generated by quadrics for all $i \geq 2$, then $I_{C_G}$ is generated by elements of the form $x_ix_j, x_i(x_j-x_k),$ 
where the variables in each product are distinct.
\end{theorem}

\begin{proof}[Proof of Theorem~\ref{thm: products}]
We proceed by induction on $d$.   The result is easy to check when $d=2.$  Assume the result for all graphs on $d-1$ vertices satisfying our hypotheses.  Our hypotheses hold for $G_i$ and $\tilde{G}_i$ for all $i \geq 2.$  Let $v= v_{i+1}.$ We may assume that $G = G_{i+1}$ and $G_i = \Gv.$

\textbf{Case 1:} \textit{$v$ has degree 1 in $G$.}
The vertex $v$ is incident to exactly one vertex $u \in \Gv$ with  $u \neq v.$  
We may assume  that  $L_v$ is spanned by a point  $p$ in $C_{\Gv}$ and the point $[0:\cdots:0:1]$. 
(By Assumption \ref{ass:label} (3), all loops are labeled by monomials.)  Then $I_{C_{\Gv}} = Q +\langle x_{d-g} \rangle,$ where $Q$ is generated by elements of the form $x_ix_j$ and $x_i(x_j-x_k)$ in which no term is divisible by $x_{d-g}$. 

We argue that $Q\subset I_{L_v}.$   Let $q = fh$ be one of the generators of $Q$ fixed above where $f$ and $h$ are linear forms.  Since $q$ must vanish at $p$, without loss of generality, we may assume that $f$ vanishes at $p.$  Since $x_{d-g}$ does not appear in $f$, then $f$ must also vanish on $[0:\cdots :1]$.   Thus, $f$ is a linear form vanishing at two points of $L_v$; hence it must vanish on all of $L_v.$  Therefore $Q \subset I_{L_v}.$  Thus, we see that $I_{C_G} = Q + \langle x_d \rangle \cdot I_{L_v}$.  Moreover, we see that $I_{C_G}$ is generated by the generators of $Q$ and elements of the form $x_{d-g}x_i$ and $x_{d-g}(x_j-x_k).$

\textbf{Case 2:} \textit{ $v$ has degree 2 in $G$.}  By Assumption \ref{ass:label} (3), there cannot be a loop at $v.$
Then there are two cases: without loss of generality, either the labels on the edges incident to $v$ have the form $x_0$ and $x_1$ or they have the form $x_0$ and $x_1-x_2.$  


In the first case, we claim that $x_0x_1 \in I_{C_{\Gv}}$.  Indeed, we have the configuration below.  

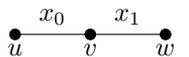
\[  \begin{tikzpicture}
	\filldraw [black] (0,0) circle (2pt);
	\filldraw [black] (1,0) circle (2pt);
	\filldraw [black] (2,0) circle (2pt);

	\draw (0,0)node[below]{$u$} -- (.5,0)node[above]{$x_0$} -- (1,0)node[below]{$v$} -- (1.5,0)node[above]{$x_1$} -- (2,0)node[below]{$w$};
		\end{tikzpicture}
		.\]
If $z$ is a vertex in $\Gv$ such that $x_0$ does not vanish on $L_z$, then $x_0$ must appear in a label on an edge incident at $z.$  If $z=u$, then $x_1$ cannot appear in a binomial on any edge incident at $z$ via Assumption \ref{ass:label} (2), and so $x_1$ vanishes on $L_z$.  If $z \neq u$, then $x_0-x_j$ must appear on an edge incident at $z.$  Again, if $x_1$ does not vanish on $L_z$, then it must appear on an edge incident at $z.$  It cannot appear in a binomial by Assumption \ref{ass:label} (1), in which case $z$ must be equal to $w$, which creates a triangle.  We conclude that either $x_0$ or $x_1$ vanishes on every irreducible component in $C_{\Gv},$ and hence that $x_0x_1 \in I_{C_{\Gv}}$.

Define a binomial minimal generator of $I_{C_{\Gv}}$ to be a binomial quadric in the ideal such that neither of its monomials is in $I_{C_{\Gv}}$.  If  $x_0(x_1-x_i)$ is in $I_{C_{\Gv}}$, then so is $x_0x_i$. Hence we may assume that we have no minimal binomial generators of the form $x_0(x_1-x_i).$  Similarly, we may assume that we have no generators of the form $x_1(x_0-x_i).$

The ideal of $I_{C_G}$ is the intersection of $I_{C_{\Gv}}$ with $I_v = \langle x_2, \ldots, x_{d-g} \rangle$, and it is generated by quadrics.  The only monomial quadrics not contained in $\langle x_2, \ldots, x_{d-g} \rangle$ are $x_0^2, x_1^2, x_0x_1.$  The monomials $x_0^2, x_1^2$ do not appear in any minimal generator of $I_{C_{\Gv}}$. Since $x_0(x_1-x_i)$ and $x_1(x_0-x_i)$ are not generators of $I_{C_{\Gv}},$ every generator of the form $x_ix_j$ and $x_i(x_j-x_k)$ must be in $I_v$ except for $x_0x_1.$  Therefore, since $I_{C_G}$ is generated by a space of quadrics whose dimension must be less than the dimension of the space of quadrics generating $I_{C_{\Gv}}$, we conclude that $I_{C_G}$ is generated by the generators of $I_{C_{\Gv}}$ minus $x_0x_1$.

In the second case, $x_0x_1$ and $x_0x_2$ are in $I_{C_{\Gv}}$ but not $I_{C_G}$ by Lemma \ref{gens} and $I_{C_G}$ is the intersection of $I_{C_{\Gv}}$ with $I_v = \langle x_1-x_2, \ldots, x_{d-g} \rangle$. We can find generators of  $I_{C_{\Gv}}$ that have the form $x_ix_j$ and $x_i(x_j-x_k).$  Note that all square-free monomials $x_ix_j$ are in $I_v$ except for $x_0x_1, x_0x_2,$ and $x_1x_2.$  The monomial $x_1x_2$ cannot be in the ideal of $C_{\Gv}$ because it contains the line parameterized by $x_1$ and $x_2.$

We claim that all of the binomial minimal generators of $I_{C_{\Gv}}$ are contained in $I_v.$  If $x_i(x_j-x_k)$ is not contained in $I_v,$ then $i$ must be 0, 1, or 2.  If it is 0, then exactly one of $j$ and $k$ is in the set $\{1,2\}.$  But, then $x_0x_1$ and $x_0x_2$ are already in $I_{C_{\Gv}}$, so $x_0(x_j-x_k)$ is not a binomial minimal generator.   

So, without loss of generality, assume $i=1$.  Let $w$ be the trivalent vertex with $L_w$ parameterized by $x_1$ and $x_2,$ and note that $w \in \Gv.$  If neither of $j$ or $k$ is in the set $\{0,2\}$ then $x_i(x_j-x_k)$ is in $I_v.$  If one of them is equal to 0, then $x_1(x_0-x_k)$ is not a binomial minimal generator because $x_0x_1 \in I_{C_{\Gv}}$.  So assume that we we have $x_1(x_2-x_k)$ with $k \neq 0, 1, 2.$  Then $x_1x_k$ is in the ideal $L_w$.  If $x_1(x_2-x_k)$ were in $I_{C_{\Gv}}$ it would also have to be in the ideal of $L_w$.  But $x_1x_k, x_1(x_2-x_k)$ in the ideal of $L_w$ would imply that  $x_1x_2$ would also be in the ideal of $L_w$, which is a contradiction.  Therefore, we have no generators of the form $x_1(x_2-x_k).$  

We conclude that all of the monomial and binomial minimal generators of $I_{C_{\Gv}}$ are in $I_{C_G}$ except for $x_0x_1$ and $x_0x_2$.  But, $x_0(x_1-x_2) \in I_{C_G},$ and the conclusion follows as in the first case since we have identified a space of quadrics in $I_{C_G}$ of dimension exactly one less than the dimension of the space of quadrics in $I_{C_{\Gv}}$.
\end{proof}

Via Corollary \ref{gleq2}, if $g \leq 2,$ we know that  each $I_{C_{G_i}}$ is generated by quadrics and we obtain the following Corollary.

\begin{corollary}\label{thm: products}
Suppose that $G$ satisfies Assumption \ref{recursive} for $p = 1$.  Fix a $\tilde{G}$ and a labeling that gives an embedding of $C_G$ into $\PP^{d-g}$ as a line arrangement.  If Assumption \ref{ass:label} is satisfied, and $g\leq 2,$ then $I_{C_G}$ is generated by elements of the form $x_ix_j, x_i(x_j-x_k),$ 
where the variables in each product are distinct.
\end{corollary}

The result in Theorem~\ref{thm: products} is sharp, as witnessed by the example below.

\begin{example}
Let $G$ be the graph below where $d=6$ and $g=2.$  Both of the vertices on the left fail part (1) of Assumption \ref{ass:label}.

\begin{figure}[htbp]
\begin{tikzpicture}
	\filldraw [black] (0,0) circle (2pt);
	\filldraw [black] (1,0) circle (2pt);
	\filldraw [black] (-1,0) circle (2pt);
	\filldraw [black] (1,1) circle (2pt);
	\filldraw [black] (-1,1) circle (2pt);
	\filldraw [black] (0,1) circle (2pt);

	\draw (-1,0) --(-.5,0)node[below]{$x_0-x_2$}--(.5,0)node[below]{$x_2$}-- (1,0) --(1,.5)node[right]{$x_4$}-- (1,1) -- (.5,1)node[above]{$x_1$}--(-.5,1)node[above]{$x_0-x_1$}--(-1,1) --(-1,.5)node[left]{$x_3$}-- cycle;
	\draw (0,0)--(0,.5)node[left]{$x_0$}--(0,1);
\end{tikzpicture}
\label{g2d6II}
\caption{}
\end{figure}
The ideal of the embedding corresponding to this labeling is
\[
 ({x}_{3} {x}_{4},{x}_{0} {x}_{4},{x}_{0} {x}_{3}-{x}_{1} {x}_{3}-{x}_{2} {x}_{3},{x}_{1} {x}_{2}).
\]
The terms in the trinomial do not appear in any of the other generators of the ideal.  Therefore, it is impossible to find a set of minimal generators that does not contain an element with at least 3 terms.
\end{example}

\begin{corollary}
If $G$ satisfies Assumption \ref{recursive} for  $p= 1$,  $g \leq 2$, and $G$ has at most 2 trivalent vertices, then there exists an embedding of $C_G \subset \PP^{d-g}$ so that $I_{C_G}$ is generated by elements of the form $x_ix_j$ and $x_i(x_j-x_k).$
\end{corollary}
\begin{proof}
Let the vertices of degree 2 with an incident edge labeled by a binomial be the last vertices in the order (so the first to get stripped of in the induction).  Combine Corollary \ref{gleq2} with Example \ref{g2} and Corollary \ref{thm: products}. 
\end{proof}


\section{Secant varieties and property $N_{k,p}$} \label{secants}

In this section we show when $N_{3,p}$ must fail for the secant line variety $\Sigma_1$.  The key idea of the proof comes from \cite{eghp}, whose authors state that their  Theorem 1.1 has a natural generalization for higher degree forms.  We give a precise statement of a special case below.

\begin{theorem}  Suppose that $X \subset \PP^n$ is a variety that satisfies $N_{k,p}.$  Let $W$ be a linear subspace of dimension $p$ with  $Z = X \cap W.$   If $\dim Z = 0,$ then $Z$ contains at most $\binom{p+k-1}{p}$ points.
\end{theorem}

\begin{proof}
It follows from the proof of Theorem 1.1 from \cite{eghp} that  the ideal of $Z $ in the homogeneous coordinate ring of $W$ is $k$-regular.  Via Theorem 4.2 in \cite{eisenbud} the degree in which the Hilbert function and Hilbert polynomial of $S_Z$ agree is the regularity of $S_Z.$  We know that the Hilbert polynomial of $S_Z$ is constant equal to the number of points in $Z.$  If $I(Z)$ is $k$-regular then $S_Z$ is $k-1$-regular.  

If $\dim (S_Z)_{k-1}$ is equal to the size of $Z,$ then $\dim S_{k-1}$ must be at least  the size of $Z.$  Hence, $|Z| \leq \binom{p+k-1}{p}.$
\end{proof}


\begin{corollary}
If $C_G$ contains a cycle of $m$ lines, then $N_{3,m-4}$ fails for the secant variety of $C_G.$
\end{corollary}
\begin{proof}
Since the $m$ lines in the cycle are contained in a $\PP^{m-1},$ so is the span of any subset of these lines.  Thus, each 3-plane obtained by taking the span of nonadjacent lines in the cycle is contained in this $\PP^{m-1}.$  There are $\binom{m}{2}-m = \frac12 m(m-3)$ such 3-planes.

A general plane of dimension $m-4$ intersects a 3-plane in $\PP^{m-1}$ in a point.  Therefore, a general $(m-4)$-plane in this $\PP^{m-1}$ intersects the secant variety of $X$ in $\frac12 m(m-3)$ points.  However, $\binom{m-4+3-1}{m-4} = \binom{m-2}{2} = \frac12 (m-2)(m-3).$   Thus, $N_{3, m-4}$ fails.
\end{proof}

\section{Questions and conjectures} \label{questions}

Computations with \emph{Macaulay 2} \cite{gs} were essential in all of our computations of embeddings of graph curves.  In addition to the results proved in this article we have several questions and conjectures regarding the defining equations and syzygies of graph curves and their secant varieties motivated by the examples that we have seen.

In \S 3 we saw that under certain hypotheses $I_{C_G}$ is generated by products of linear forms that can be described explicitly in terms of the combinatorics of the graph $G.$  The combinatorics of the $k$th secant variety of $C_G$ is encoded in an \emph{intersection lattice} whose elements are constructed by intersecting subsets of the subspaces.  From the intersection lattice of an arrangement, we get a partially ordered set ordered by reverse inclusion of subspaces. 
\begin{question}
Does the partially ordered set associated to the $k$th secant variety have any interesting combinatorial features?  We conjecture that $\Sigma_k$ is Cohen-Macaulay, so will the corresponding poset be shellable?
\end{question}

It is also natural to ask if there is an analogue of Theorem~\ref{thm: products} for secant varieties, perhaps requiring additional hypotheses on the intersection lattice of the secant varieties of $C_G.$
\begin{question}\label{question:pl}
Are the secant varieties of $C_G$ defined by products of linear forms?
\end{question}

Finding generators of $I_{\Sigma_k}$ that are products of linear forms is equivalent to finding an explicit and special basis for the ideal that may have combinatorial interest.  Of course, a module does not typically have a unique generating set or a unique minimal free graded resolution.  However, the number of minimal generators of degree $j$ of the $i$th syzygy module is invariant under a change of basis.   Given a finitely generated graded module $M$, the graded Betti number $\beta_{i,j}$ is the number of minimal generators of degree $j$ required at the $i$th stage of a minimal free graded resolution of $M.$  A standard way of displaying the graded Betti numbers of a module is with a graded Betti diagram organized as follows:

\[
\begin{array}{c|cccc}
& 0 & 1 & 2 & \\
\hline
0&\beta_{0,0} & \beta_{1,1} & \beta_{2,2}& \cdots \\
1&\beta_{0,1} & \beta_{1,2} & \beta_{2,3}&\cdots
\end{array}
\]

Bounds on the number of rows and columns of the graded Betti diagram of a module give a rough sense of how complicated it is.  Specifically, recall that the \emph{regularity} of a finitely generated graded module $M$ is equal to $\sup\{ j-i \mid \beta_{i,j} \neq  0 \text{ for some } i\},$ and thus regularity gives a bound on the number of rows of the graded Betti diagram of $M.$   Additionally, by the Auslander-Buchsbaum formula, a variety $X \subset \PP^n$ is arithmetically Cohen-Macaulay if $\sup\{ i \mid \beta_{i,j} \neq 0 \text{ for some } i\} = \codim X$, which bounds the number of columns of the graded Betti diagram off $M$.  The following conjecture is the graph curve analogue of Conjecture 1.4 in \cite{sidmanVermeire} which refines conjectures from \cite{vermeire}.
 
\begin{conjecture}\label{conj:sec}
If Assumption \ref{recursive} holds for some $p \geq 2k$, then the $k$th secant variety of $C_G$ has regularity equal to $2k+1$ and is arithmetically Cohen-Macaulay.  \end{conjecture}

Note that as the secant varieties of $C_G$ are not normal, we cannot expect projective normality.  

In addition to bounding the length and width of the graded Betti diagram we conjecture that under certain conditions, one particular graded Betti number counts the number of cycles of minimal length in the graph.  Recall that the \emph{girth} of a graph is the length of its smallest cycle.

\begin{conjecture}\label{conj:betti}
Let $G $ be a graph on $d$ vertices, embedded as in Theorem 1.3.  Let $n$ denote the girth of $G$. Assume that $d = 2g+1+p$ and $n-2 \leq p.$ Then  property $N_p$ fails and $\beta_{n-2,n}$ is equal to the number of cycles of length $n$ in $G$.
\end{conjecture}


Example~\ref{ex:g2d10} gives an illustration of the properties discussed in Conjectures~\ref{conj:sec}, \ref{conj:betti}.

\begin{example}[$g = 2,$ $d = 10$]\label{ex:g2d10}

Let $G$ be as given in Figure~\ref{fig:graph}.

\begin{figure}[htbp]
\begin{tikzpicture}[style=thick]

\node(1) at (-1.5, -1) [label=left:]{};
\node (2) at (-1.5, 0) [label=left:]{};

\node (3) at (-1.5, 1) [label=left:]{};

\node(4) at (0, 1.5) [label=above:]{};

\node(5) at (1.5, 1) [label=right:]{};

\node (6) at (1.5, 0) [label=right:]{};

\node (7) at (1.5, -1) [label=right:]{};

\node (8) at (0, -1.5) [label=below:]{};

\node(9) at (0, -.60) [label=left:]{};
\node (10) at (0, .60) [label=left:]{};

\filldraw (1) circle (4pt);
\filldraw (2) circle (4pt);
\filldraw (3) circle (4pt);
\filldraw (4) circle (4pt);
\filldraw (5) circle (4pt);
\filldraw (6) circle (4pt);
\filldraw (7) circle (4pt);
\filldraw (8) circle (4pt);
\filldraw (9) circle (4pt);
\filldraw (10) circle (4pt);

\draw (1)-- (-1.5, -.5)node[left]{$x_7$}--(2)--(-1.5, .7)node[left]{$x_8$}--(3)--(-.75, 1.25)node[above]{$x_0-x_6$}--(4)--(.75, 1.25)node[above]{$x_0$}--(5)--(1.5, .7)node[right]{$x_1$}--(6)--(1.5, -.5)node[right]{$x_2$}--(7)--(.75, -1.25)node[below]{$x_3$}--(8)--(-.75, -1.25)node[below]{$x_3-x_4$}--(1);
\draw (8) -- (0, -.8)node[left]{$x_4$}--(9)--(0, 0)node[left]{$x_5$}--(10)--(0, .9)node[left]{$x_6$}--(4);
\end{tikzpicture}
\caption{}
\label{fig:graph}
\end{figure}
The ideal of $C_G$ corresponding to this labeling is given below.
\[
I_{C_G} =  \begin{matrix}
({x}_{5} {x}_{8},&
 {x}_{4} {x}_{8},& 
 {x}_{3} {x}_{8},& 
 {x}_{2} {x}_{8},& 
 {x}_{1} {x}_{8},& 
{x}_{6} {x}_{7}, &
      {x}_{5} {x}_{7},\\
      
        {x}_{2} {x}_{7},& 
       {x}_{1} {x}_{7},& 
       {x}_{0} {x}_{7},& 
       {x}_{4} {x}_{6},&
        {x}_{3} {x}_{6},&
      {x}_{2}{x}_{6},&
       {x}_{1} {x}_{6},\\
        {x}_{3} {x}_{5},&
         {x}_{2} {x}_{5},&
          {x}_{1} {x}_{5},&
           {x}_{0} {x}_{5},&
            {x}_{2} {x}_{4},&
      {x}_{1} {x}_{4},&
      {x}_{0} {x}_{4},\\
      {x}_{1} {x}_{3},&
      {x}_{0} {x}_{3},&
      {x}_{0} {x}_{2},&
      {x}_{3} {x}_{7}-{x}_{4} {x}_{7},&

       {x}_{0} {x}_{8}-{x}_{6} {x}_{8})
 \end{matrix}
 \]
 
The graded Betti diagram of $S/I_{C_G}$ shows that $N_{2,5}$ fails as $\beta_{5,7} = 2$.  As Conjecture~\ref{conj:betti} predicts, the girth of $G$ is 7, and $G$ contains precisely 2 cycles of length 7.
 \[\begin{matrix}
      &0&1&2&3&4&5&6&7\\\text{total:}&1&26&98&168&154&72&15&2\\\text{0:}&1&\text{.}&\text{.}&\text{.}&\text{.}&\text{.}&\text{.}&\text{.}\\\text{1:}&\text{.}&26&98&168&154&70&8&\text{.}\\\text{2:}&\text{.}&\text{.}&\text{.}&\text{.}&\text{.}&2&7&2\\\end{matrix}.\]
      
 We can also compute the ideal of $\Sigma$
      
 \[I(\Sigma) = 
 \begin{matrix}(
 {x}_{3} {x}_{5} {x}_{8},&
 {x}_{2} {x}_{5} {x}_{8},&
 {x}_{1} {x}_{5} {x}_{8},&
 {x}_{0} {x}_{4} {x}_{8}-{x}_{4} {x}_{6} {x}_{8},\\
 
 {x}_{2} {x}_{4} {x}_{8},&
 {x}_{1} {x}_{4}{x}_{8},&
 {x}_{1} {x}_{3} {x}_{8},&
 {x}_{0} {x}_{3} {x}_{8}-{x}_{3} {x}_{6} {x}_{8},\\
 
 {x}_{2} {x}_{6} {x}_{7},&
 {x}_{1} {x}_{6} {x}_{7},&
 {x}_{2} {x}_{5}{x}_{7},&
 {x}_{0} {x}_{2}{x}_{8}-{x}_{2} {x}_{6} {x}_{8},\\
 
 {x}_{1} {x}_{5} {x}_{7},&
 {x}_{0} {x}_{5} {x}_{7},&
 {x}_{0} {x}_{2} {x}_{7},&
 {x}_{3} {x}_{6} {x}_{7}-{x}_{4} {x}_{6} {x}_{7},\\

 {x}_{2} {x}_{4} {x}_{6},&
 {x}_{1} {x}_{4} {x}_{6},&
 {x}_{1} {x}_{3} {x}_{6},&
 {x}_{1} {x}_{3} {x}_{7}-{x}_{1} {x}_{4} {x}_{7},\\

 {x}_{1} {x}_{3} {x}_{5},&
 {x}_{0} {x}_{3}{x}_{5},&
 {x}_{0} {x}_{2} {x}_{5},&
 {x}_{0} {x}_{3} {x}_{7}-{x}_{0} {x}_{4}{x}_{7},\\

 {x}_{0} {x}_{2} {x}_{4}) \end{matrix}
 \]
 and its graded Betti diagram
   \[\begin{matrix}
      &0&1&2&3&4&5\\\text{total:}&1&25&58&43&12&3\\\text{0:}&1&\text{.}&\text{.}&\text{.}&\text{.}&\text{.}\\\text{1:}&\text{.}&\text{.}&\text{.}&\text
      {.}&\text{.}&\text{.}\\\text{2:}&\text{.}&25&58&41&\text{.}&\text{.}\\\text{3:}&\text{.}&\text{.}&\text{.}&\text{.}&7&\text{.}\\\text{4:}&\text{.}&\text{.}&\text{.}&2&5&3\\\end{matrix}.\]
We see that $N_{3,3}$ fails for $\Sigma$ and that $\beta_{3,7} =2,$ which is the number of cycles of length equal to the girth of $G.$  We can also see from the graded Betti diagram that $\Sigma$ is arithmetically Cohen-Macaulay and that $I(\Sigma)$ has regularity 5.
 \end{example}

It is natural to ask if combinatorics can be used to compute other values of the $\beta_{i,j}.$  One result that gives the flavor of what might be possible is due to Gasharov, Peeva and Welker \cite{gpw} who used the lcm lattice of a monomial ideal to compute graded Betti numbers of monomial ideals.  

\begin{question}Is there an analogue of the lcm lattice for graph curves and their secant varieties that would allow us to compute (or estimate) the graded Betti numbers of graph curves?
\end{question}

Further work on understanding the graded Betti numbers of graph curves has been done by \cite{bknpv}.

It is also interesting to consider $C_G$ as a deformation of a smooth curve.  In Example~\ref{ex:g2d10}, $C_G$ has a $7$-secant $\PP^5$ while a smooth curve of genus $2$ in $\PP^8$ has no such $\PP^5$.  As any strictly subtrivalent graph curve $C_G\subset\PP^n$ is smoothable in $\PP^n$ \cite[29.9]{hartdef}, it is our expectation that we have a family of seven $6$-secant $\PP^5$s to smooth curves that collapse to the $7$-secant $\PP^5$ in the singular limit $C_G$.  It also seems reasonable to believe that the secant varieties to embedded curves in a flat family themselves form a flat family, and so the secant varieties to $C_G$ should, in particular, have the same dimension and degree as those to smooth curves.  In fact, since each pair of disjoint lines in $C_G$ spans a $\PP^3,$ we have a 3-dimensional secant plane for each edge in the complement of the graph $G.$  If $C_G$ has degree $d$ and genus $g$, then $G$ has $d$ vertices and $d+g-1$ edges.  Thus, the number of edges in the complement of $G$ is $\binom{d}{2} - d-g+1 = \binom{d-1}{2}-g,$  which is the degree of the secant variety of a smooth curve of degree $d$ and genus $g.$


\bibliographystyle{plain}
\bibliography{brsvArxiv}

\end{document}